\documentclass[a4paper,12pt]{article}
\usepackage{amssymb}
\usepackage{amscd}

\title{Wild hypersurface bundles over toric varieties}

\date{}

\author{{\sc Hiroshi Sato}\thanks{Partly supported by the 
Grant-in-Aid for JSPS Fellows, The Ministry of Education, 
Science, Sports and Culture, Japan.
\newline
\hspace*{1.5em} {\em $2000$ Mathematics Subject Classification\/}.
Primary 14M25;
Secondary 14E30, 14J40, 14J45.}}

\newtheorem{Thm}{Theorem}[section]
\newtheorem{Prop}[Thm]{Proposition}
\newtheorem{Cor}[Thm]{Corollary}

\newtheorem{Def}[Thm]{Definition}
\newtheorem{Rem}[Thm]{Remark}
\newtheorem{Ex}[Thm]{Example}

\newcommand{\proof}{Proof. \quad}
\newcommand{\qed}{\hfill q.e.d.}

\newcommand{\Pic}{\mathop{\rm Pic}\nolimits}

\newcommand{\G}{\mathop{\rm G}\nolimits}

\newcommand{\PC}{\mathop{\rm PC}\nolimits}

\begin{document}
\maketitle

\begin{abstract}                
In this paper, we investigate when there exists 
a wild hypersurface bundle over a smooth proper 
toric variety in positive characteristic. 
In particular, we determine the possibilities for 
toric varieties with Picard number at most three or 
toric Fano varieties of dimension at most four. 
Moreover, we can construct wild hypersurface 
bundles over them.

\end{abstract}

\thispagestyle{empty}

\section{Introduction}

\hspace{5mm} A {\em wild hypersurface bundle} is 
a peculiar phenomenon in positive characteristic 
(see Definition \ref{defwhb}). 
Only few examples of wild hypersurface bundles 
are known. Saito \cite{saito1} 
completely determined when a smooth 
Fano $3$-folds with Picard 
number $2$ has a wild conic bundle structre. 
As a generalization for this result, Mori-Saito 
\cite{mori1} showed the following:

\begin{Thm}[Mori-Saito \cite{mori1}]
Let $f:X\to S$ be a wild 
hypersurface bundle of degree $p$, 
$d=\dim S$ and $\dim X=2d-1$. 
If $S$ is isomorphic to 
a direct product of projective spaces, 
then one of the following holds$:$
\begin{enumerate}
\item $S\simeq{\mathbb P}^d$ and $X$ is a 
smooth divisor of bidegree $(1,p)$ in 
${\mathbb P}^d\times{\mathbb P}^d$.
\item $p=2$, $S\simeq({\mathbb P}^1)^d$ and 
$X$ is a smooth divisor in 
$Y={\mathbb P}_S(\mathcal{O}_S\oplus
\bigoplus_{i=1}^{d}p_i^{*}
\mathcal{O}_{{\mathbb P}^1}(1))$ such that 
$X\sim 2\xi$, where $p_i:S\to{\mathbb P}^1$ is the 
$i$-th projection and $\xi$ is the tautological 
line bundle of $Y\to S$.
\end{enumerate}
\end{Thm}

In this paper, we consider the case where $S$ is a 
smooth proper toric $d$-fold. 
Using the technique in Mori-Saito \cite{mori1}, we completely 
determine the possibilities for $S$ when the Picard 
number of $S$ is $2$ or $3$ (see Section \ref{pic2or3}), 
or $S$ is a toric Fano $d$-fold with $d\leq 4$ 
(see Section \ref{overfano}). Moreover, we can construct 
wild hypersurface bundles for these cases. 

The content of this paper is as follows: Section \ref{junbi} is 
a section for preparation. We review the concepts of primitive 
collections and relations, and explicitly describe the fans 
for projective space bundles over toric varieties. In Section 
\ref{whbsec}, we review the definition of wild hypersurface 
bundles. The combinatorial version of the key result in 
Mori-Saito \cite{mori1} is given. In Section \ref{pic2or3}, 
we consider the case where the Picard number of $S$ is $2$ or 
$3$. There exist two new classes which have 
wild hypersurface bundle structures. In Section \ref{overfano}, 
we consider the case where $S$ is a toric Fano variety. 
In particular, we determine the toric Fano $d$-folds which 
have wild hypersurface bundle structures for $d\leq 4$. 
These Fano varieties are interesting from the viewpoint of 
the birational geometry. 

The author would like to thank Doctor Natsuo Saito for 
introducing the author to this problem 
and giving useful comments. 
The author also wishes to thank Professors 
Shihoko Ishii and Osamu Fujino 
for advice and encouragement. 

\section{Preliminaries}\label{junbi}

\hspace{5mm} This section is devoted to explaining some basic 
facts of the toric geometry. See Batyrev \cite{batyrev3}, 
\cite{batyrev4}, Fulton \cite{fulton1}, Oda \cite{oda2} and 
Sato \cite{sato1} more precisely.

Let $S=S_{\Sigma}$ be a smooth proper toric $d$-fold 
associated 
to a fan $\Sigma$ over an algebraically closed field. 
Let $\G(\Sigma)$ be the set of primitive generators of 
$1$-dimensional cones in $\Sigma$. A subset $P\subset 
\G(\Sigma)$ is called a {\em primitive collection} 
if $P$ does not generate a cone in $\Sigma$, while 
any proper subset of $P$ generates a cone in $\Sigma$. 
We denote by $\PC(\Sigma)$ the set of primitive 
collections of $\Sigma$. 
For a primitive collection $P=\{x_1,\ldots,x_m\}$, there exists 
the unique cone $\sigma(P)$ in $\Sigma$ such that $x_1+\cdots+x_m$ 
is contained in the relative interior of $\sigma(P)$, 
since $S$ is proper. So, we obtain an equality
$$x_1+\cdots+x_m=b_1y_1+\cdots+b_ny_n,$$
where $y_1,\ldots,y_n$ are the generators of $\sigma(P)$, 
that is, $\sigma(P)\cap\G(\Sigma)=\{y_1,\ldots,y_n\}$, 
and $b_1,\ldots,b_n$ are positive integers. 
We call this equality the {\em primitive relation} of $P$. 
Thus, we obtain an element $r(P)$ in $A_{1}(S)$ for any 
primitive collection $P\in\PC(\Sigma)$, where $A_{1}(S)$ 
is the group of $1$-cycles on $S$ modulo rational 
equivalences. We define the {\em degree} of $P$ as 
$\deg P:=\left( -K_S\cdot r(P)\right)=m-(a_1+\cdots+a_n)$. 

\begin{Prop}[Batyrev \cite{batyrev3}, 
Reid \cite{reid1}]\label{toriccone}

Let $S=S_{\Sigma}$ be a smooth projective toric variety. Then
$${\rm NE}(S)=\sum_{P\in\PC(\Sigma)}{\mathbb R}_{\geq 0}r(P),$$
where ${\rm NE}(S)$ is the Mori cone of $S$.

\end{Prop}

A primitive collection $P$ is said to be 
{\em extremal} if $r(P)$ is contained 
in an extremal ray of ${\rm NE}(S)$. For the torus 
invariant curve $C$ 
contained in this extremal ray, we have 
$$N_{C/S}\simeq \mathcal{O}_C(1)^{\oplus(m-2)}
\oplus\mathcal{O}_C^{\oplus(d-m-n+1)}
\oplus\mathcal{O}_C(-b_1)\oplus\cdots\oplus\mathcal{O}_C(-b_n),$$
where $N_{C/S}$ is the normal bundle.

\bigskip

Next, we explain how to construct the fan
corresponding to a projective space bundle over 
a toric variety.

Let $S=S_{\Sigma}$ be a smooth proper toric $d$-fold, 
$\Sigma$ a fan in $N={\mathbb Z}^d$, 
$\G(\Sigma)=\{x_1,\ldots,x_l\}$ and $D_1,\ldots,D_l$ 
the torus invariant prime divisors corresponding to 
$x_1,\ldots,x_l$, respectively. 
For torus invariant divisors 
$$E_1=\sum_{i=1}^lc_{1,i}D_i,\ \ldots,
\ E_r=\sum_{i=1}^lc_{r,i}D_i,$$
put 
$$E=\mathcal{O}\oplus
\mathcal{O}_S(E_1)\oplus\cdots\oplus\mathcal{O}_S(E_r).$$
We construct the fan $\widetilde{\Sigma}$ in 
$\widetilde{N}:=N\oplus {\mathbb Z}^r$ 
corresponding to the ${\mathbb P}^r$-bundle 
${\mathbb P}_S(E)$ over $S$.

Let $\{e_1,\ldots,e_r\}$ be the standard basis for 
${\mathbb Z}^r$. The elements of $\G(\widetilde{\Sigma})$ are 
$$y_1:=e_1,\ \ldots,\ y_r:=
e_r,\ y_{r+1}:=-(e_1+\cdots+e_r),\ $$
$$\widetilde{x}_1:=x_1+\sum_{i=1}^rc_{i,1}e_i,\ \ldots,
\ \widetilde{x}_l:=x_l+\sum_{i=1}^rc_{i,l}e_i.$$
For a maximal cone $\sigma={\mathbb R}_{\geq 0}x_{i_1}+\cdots
+{\mathbb R}_{\geq 0}x_{i_d}$ in $\Sigma$, put 
$\widetilde{\sigma}:={\mathbb R}_{\geq 0}\widetilde{x}_{i_1}+
\cdots+{\mathbb R}_{\geq 0}\widetilde{x}_{i_d}
\subset\widetilde{N}\otimes{\mathbb R}$. 
Put $\widetilde{\tau}_i:=
{\mathbb R}_{\geq 0}y_1+\cdots+
{\mathbb R}_{\geq 0}y_{i-1}+
{\mathbb R}_{\geq 0}y_{i+1}+\cdots+
{\mathbb R}_{\geq 0}y_{r+1}
\subset\widetilde{N}\otimes{\mathbb R}$ for $1\leq i\leq r+1$. 
The set of maximal cones in $\widetilde{\Sigma}$ is
$$\left\{\left.\widetilde{\sigma}+\widetilde{\tau}_i
\,\right|\,\sigma\ \mbox{is a maximal cone in}\ \Sigma,\ 
1\leq i\leq r+1\right\}.$$
The tautological line bundle $\xi$ for ${\mathbb P}_S(E)\to S$ 
is $\mathcal{O}_{{\mathbb P}_S(E)}(F_{r+1})$, 
where $F_{r+1}$ is the torus invariant prime divisor 
corresponding to $y_{r+1}$.

\section{Wild hypersurface bundles}\label{whbsec}

\hspace{5mm} In this section, we review the definition 
of a {\em wild hypersurface bundle structure} 
and some results in Mori-Saito \cite{mori1}. 
From now on, we work over an algebraically 
closed field $k$ of characteristic $p>0$.

\begin{Def}[Mori-Saito \cite{mori1}]\label{defwhb}

{\rm
Let $X$ and $S$ be smooth algebraic varieties over $k$, 
and $f:X\to S$ a projective flat morphism with $M$  a 
relatively very ample divisor such that $X$ is embedded in 
$\pi:{\mathbb P}_S(E)\to S$, where $E=f_{*}M$. We call $f$ a 
{\em wild hypersurface bundle of degree $p$} if for any 
$s\in S$, the geometric fiber $f^{-1}(s)$ is defined in 
${\mathbb P}_S(E)$ by $x^p=0$ for some non-zero $x\in E_s$.
}

\end{Def}

Let $\xi$ be the tautological line bundle of 
${\mathbb P}_S(E)$. Then, there exists a Cartier divisor 
$L$ on $S$ such that $X\sim p\xi+\pi^{*}L$ in 
$\Pic{\mathbb P}_S(E)$. 
Let $d=\dim S$. If $\dim X=2d-1$, then there exists 
an exact sequence 
\begin{equation}
0\to\mathcal{O}_S\to E^p\otimes L\to T_S\to 0,
\end{equation}
where $T_S$ is the tangent bundle of $S$ 
(see Theorem 1 in Mori-Saito \cite{mori1}). 
Thanks to this exact sequence, we can study 
wild hypersurface bundle structures easily. 
So, in this paper, we add the assumption $\dim X=2d-1$ 
to the definition of a wild hypersurface bundle of 
degree $p$. 
We will use these notation throughout this paper.

The following is a slight generalization 
of Proposition $5$ in Mori-Saito \cite{mori1}. 
The proof is similar.

\begin{Prop}\label{keyprop}

Let $f:X\to S$ be a wild hypersurface bundle of degree $p$ 
and $C$ a normal rational curve on $S$ such that
$$T_S\otimes \mathcal{O}_C\simeq \bigoplus_{i=-\infty}^{2}
\mathcal{O}_C(i)^{\oplus a_i}.$$
Then, the following hold.
\begin{enumerate}
\item If the restriction of the exact sequence $(1)$ on $C$ is 
non-split, then for any $a_i>0$, $i-1$ is divisible by 
$p$.
\item If the restriction of the exact sequence $(1)$ on $C$ is 
split, then $p=2$ and for any $a_i>0$, $i$ is 
an even number.
\end{enumerate}

\end{Prop}

\begin{Rem}\label{normbdl}

{\rm
In Proposition \ref{keyprop},
$$E^p\otimes L\otimes\mathcal{O}_C\simeq
\left(\bigoplus_{i=-\infty}^{-1}\mathcal{O}_C(i)^{\oplus a_i}
\right)\oplus\mathcal{O}_C^{\oplus a_0}\oplus
\mathcal{O}_C(1)^{\oplus(a_1+2)}\oplus
\mathcal{O}_C(2)^{\oplus(a_2-1)}$$
for the case (i), while
$$E^p\otimes L\otimes\mathcal{O}_C\simeq
\left(\bigoplus_{i=-\infty}^{-1}\mathcal{O}_C(i)^{\oplus a_i}
\right)\oplus\mathcal{O}_C^{\oplus(a_0+1)}\oplus
\mathcal{O}_C(1)^{\oplus a_1}\oplus
\mathcal{O}_C(2)^{\oplus a_2}$$
for the case (ii).
}

\end{Rem}

We apply this result for the case where $S$ is a toric variety.

\begin{Cor}\label{torickey}

Let $S=S_{\Sigma}$ be a smooth proper 
toric $d$-fold and $f:X\to S$ 
a wild hypersurface bundle of degree $p$. For an extremal 
primitive relation 
$$x_1+\cdots+x_m=b_1y_1+\cdots+b_ny_n,$$
where $\{x_1,\ldots,x_m,y_1,\ldots,y_n\}\subset\G(\Sigma)$ and 
$b_1,\ldots,b_n$ are positive integers, one of the following 
holds.
\begin{enumerate}
\item $m+n=d+1$ and $b_i+1$ is divisible by $p$ for any $i$. 
\item $p=2$, $m=2$ and $b_i$ is an even number 
for any $i$. 
\end{enumerate}

\end{Cor}
\proof
This can be proven by Proposition \ref{keyprop} immediately. 
For the case (i), $N_{C/S}$ does not contained 
$\mathcal{O}_C$, so $m+n=d+1$. The left part is similar.
\qed

\begin{Rem}
{\rm
For the case (i) in Corollary \ref{torickey}, let 
$\varphi:S\to\overline{S}$ be the associated 
extremal contraction. If $S\not\simeq\mathbb{P}^d$, then 
$\varphi$ is birational and 
the image of the exceptional set of $\varphi$ is 
a point. 
}
\end{Rem}

\section{Toric varieties with Picard number 2 or 3}
\label{pic2or3}

\hspace{5mm} In this section, we treat the case where $S$ 
is a smooth proper toric $d$-fold with Picard number 
$2$ or $3$. We construct some examples of wild 
hypersurface bundles using the 
notion of {\em homogeneous coordinate rings} of 
toric varieties (see Cox \cite{cox1}). 

\bigskip

(I) The Picard number of $S$ is two.

\begin{Prop}\label{toricpic2}
Let $S$ be a smooth proper toric $d$-fold with 
Picard number $2$. If there exists a wild hypersurface 
bundle $f:X\to S$, then $p=2$ and $S$ is isomorphic to 
either ${\mathbb P}^1\times{\mathbb P}^1$ or 
$${\mathbb P}_{{\mathbb P}^{d-1}}
(\mathcal{O}_{{\mathbb P}^{d-1}}\oplus
\mathcal{O}_{{\mathbb P}^{d-1}}(2a-1)),$$
where $a$ is a positive integer.
\end{Prop}
\proof
There exists a wild hypersurface bundle of degree $2$ over 
${\mathbb P}^1\times{\mathbb P}^1$ 
(see Mori-Saito \cite{mori1}). So, suppose 
$S\not\simeq{\mathbb P}^1\times{\mathbb P}^1$. 

$S$ is a ${\mathbb P}^m$-bundle over ${\mathbb P}^n$ by 
the classification of proper toric varieties with 
Picard number $2$ (see Kleinschmidt \cite{klein1}). 
On the other hand, $m=1$ by the case (ii) in 
Proposition \ref{keyprop}. Thus, $p=2$ and $S\simeq
{\mathbb P}_{{\mathbb P}^{d-1}}
(\mathcal{O}_{{\mathbb P}^{d-1}}\oplus
\mathcal{O}_{{\mathbb P}^{d-1}}(\alpha))$ 
for a non-negative integer $\alpha$. 
Since the normal bundle $N_{C_1/S}$ of 
the torus invariant curve $C_1$ contained in another extremal 
ray is $$\mathcal{O}_{C_{1}}(-\alpha)\oplus
\mathcal{O}_{C_{1}}(1)^{\oplus(d-2)},$$
$\alpha$ is an odd number by the case (i) in Proposition 
\ref{keyprop}.
\qed

\bigskip

Next, we construct a wild hypersurface bundle 
of degree $2$ for the above case. So, let 
$S:={\mathbb P}_{{\mathbb P}^{d-1}}
(\mathcal{O}_{{\mathbb P}^{d-1}}\oplus
\mathcal{O}_{{\mathbb P}^{d-1}}(2a-1))$ for a 
positive integer $a$ and $\Sigma$ the associated fan. 
Then, the primitive relations of $\Sigma$ are 
$$({\rm a})\ x_1+\cdots+x_d=(2a-1)x_{d+1}\ \mbox{and}
\ ({\rm b})\ x_{d+1}+x_{d+2}=0,$$
where $\G(\Sigma)=\{x_1,\ldots,x_{d+2}\}$. Let $D_1,\ldots,
D_{d+2}$ be the torus invariant prime divisors corresponding 
to $x_1,\ldots,x_{d+2}$, respectively. 
We may assume that $\{x_1,\ldots,x_{d-1},x_{d+1}\}$ is
the standard basis for $N$. 
By considering the divisors of the rational functions 
corresponding to $x_1,\ldots,x_{d-1},x_{d+1}$, we have 
$D_1=\cdots=D_d$ and $D_{d+2}=(2a-1)D_1+D_{d+1}$ in
 $\Pic S$. 
Let $C_1$ and $C_2$ be the torus invariant curves corresponding 
to the extremal primitive relations (a) and (b), respectively. 
Then, $(D_1\cdot C_1)=1$, $(D_{d+1}\cdot C_1)=-(2a-1)$, 
$(D_1\cdot C_2)=0$ and $(D_{d+1}\cdot C_2)=1$. 
Put 
$$E=\mathcal{O}_S^{\oplus d}\oplus
\mathcal{O}_S((a-1)D_1+D_{d+1})\ \mbox{and}\ 
L=\mathcal{O}_S(D_1).$$
Then, we can easily check that $E$ and $L$ satisfy 
the conditions 
$$E^2\otimes L\otimes\mathcal{O}_{C_1}=
\mathcal{O}_{C_1}(-1)\oplus\mathcal{O}_{C_1}(1)^{\oplus d}\ 
\mbox{and}\ E^2\otimes L\otimes\mathcal{O}_{C_2}=
\mathcal{O}_{C_2}^{\oplus d}\oplus\mathcal{O}_{C_2}(2).$$
In fact, we can construct a wild hypersurface bundle 
for these $E$ and $L$ as follows.

Let $\widetilde{\Sigma}$ be the fan corresponding to 
$Y=\mathbb{P}_S(E)$. We use the same notation as in 
Section \ref{junbi}. The primitive relations of 
$\widetilde{\Sigma}$ are $\widetilde{x}_{d+1}+
\widetilde{x}_{d+2}=y_{1}$, 
$y_{1}+\cdots+y_{d+1}=0$ and
$$\widetilde{x}_{1}+\cdots+\widetilde{x}_{d}
=\left\{ \begin{array}{ccl}
(2a-1)\widetilde{x}_{d}+y_{2}+\cdots
+y_{d+1} & 
\mbox{if} & a=1 \\
(2a-1)\widetilde{x}_{d}+(a-2)y_{1} & 
 & \mbox{otherwise}, \\
\end{array} \right. $$
where $\G(\widetilde{\Sigma})=\{\widetilde{x}_1,\ldots,
\widetilde{x}_{d+2},y_1,\ldots,
y_{d+1}\}$. Let $\widetilde{D}_1,\ldots,
\widetilde{D}_{d+2},F_1,\ldots,
F_{d+1}$ be the torus invariant prime 
divisors corresponding to $\widetilde{x}_1,\ldots,
\widetilde{x}_{d+2},y_1,\ldots,
y_{d+1}$, respectively. Then, we have 
$\widetilde{D}_1=\cdots=\widetilde{D}_{d}$, 
$\widetilde{D}_{d+2}=(2a-1)\widetilde{D}_d+
\widetilde{D}_{d+1}$, $F_2=\cdots=
F_{d+1}$ and 
$F_{d+1}=(a-1)\widetilde{D}_1+\widetilde{D}_{d+1}
+F_1$ in $\Pic Y$. 
Since the tautological line bundle $\xi$ for $\pi:Y\to S$ 
is $\mathcal{O}_Y(F_{r+1})$, we have 
$X\sim 2\xi+\pi^{*}L=2F_{d+1}+\widetilde{D}_1
=\widetilde{D}_{d+1}+\widetilde{D}_{d+2}+
2F_{1}$. 
Thus, for example, the smooth 
hypersurface $X$ in $Y$ defined by the equation
$$X_{d+1}X_{d+2}Y_1^2+X_1Y_2^2+\cdots+X_dY_{d+1}^2=0$$
is a wild hypersurface 
bundle of degree $2$ over $S$, 
where $X_1,\ldots,X_{d+2},Y_1,\ldots,Y_{d+1}$ are 
the homogeneours coordinates of $Y$ corresponding to
$\widetilde{D}_1,\ldots,
\widetilde{D}_{d+2},F_1,\ldots,
F_{d+1}$, respectively. We can easily check the 
smoothness of $X$, so we leave the details for the exercise.

\bigskip

(II) The Picard number of $S$ is three.

\bigskip

In this case, we suppose $d\geq 3$.

Batyrev \cite{batyrev3} classified smooth projective 
toric $d$-folds with Picard number $3$ 
using the notion of primitive relations. 

\begin{Thm}[Batyrev \cite{batyrev3}]\label{tmjtmj}

Let $S=S_{\Sigma}$ be a smooth projective toric $d$-fold with 
Picard number three. Then, one of the following holds.
\begin{enumerate}
\item $\#PC(\Sigma)=3$, and for any distinct elements 
$P_1,P_2\in\PC(\Sigma)$, we have $P_1\cap P_2=\emptyset.$
\item $\#\PC(\Sigma)=5$, and there exists 
$(p_0,p_1,p_2,p_3,p_4)\in({\mathbb Z}_{>0})^{5}$ 
such that $p_0+p_1+p_2+p_3+p_4=d+3$ and 
the primitive relations of $\Sigma$ are
$$v_1+\cdots+v_{p_0}+y_{1}+\cdots+y_{p_1}=
c_2z_2+\cdots+c_{p_2}z_{p_2}+(b_1+1)t_1+
\cdots+(b_{p_3}+1)t_{p_3},$$
$$y_1+\cdots+y_{p_1}+z_1+\cdots+z_{p_2}=
u_1+\cdots+u_{p_4},\ z_1+\cdots+z_{p_2}+t_1+\cdots+t_{p_3}=0,$$
$$t_1+\cdots+t_{p_3}+u_1+\cdots+u_{p_4}=
y_{1}+\cdots+y_{p_1}\mbox{ and}$$
$$u_1+\cdots+u_{p_4}+v_1+\cdots+v_{p_0}=
c_2z_2+\cdots+c_{p_2}z_{p_2}+b_1t_1+\cdots+b_{p_3}t_{p_3},$$
where
$$\G(\Sigma)=\{v_1,\ldots,v_{p_0},y_{1},
\ldots,y_{p_1},z_1,\ldots,z_{p_2},t_1,\ldots,
t_{p_3},u_1,\cdots,u_{p_4}\}$$
and $c_2,\ldots,c_{p_2},b_1,\ldots,b_{p_3}$ 
are positive integers.
\end{enumerate}
\end{Thm}

For positive integers $a$ and $b$, 
let $\Sigma^{d}(a,b)$ be the fan whose 
primitive relations are 
$$x_1+\cdots+x_{d-1}=(2a-1)x_d+(2b-1)x_{d+2},\ x_d+x_{d+1}=0
\ \mbox{and}\ x_{d+2}+x_{d+3}=0,$$
where $\G(\Sigma^{d}(a,b))=\{x_1,\ldots,x_{d+3}\}$ and 
$W^d(a,b)$ the associated toric $d$-fold with 
Picard number $3$. The following Proposition holds.

\begin{Prop}\label{toricpic3}
Let $S$ be a smooth proper toric $d$-fold with 
Picard number $3$. If there exists a wild hypersurface 
bundle $f:X\to S$, then $p=2$ and $S$ is isomorphic to 
either ${\mathbb P}^1\times{\mathbb P}^1\times{\mathbb P}^1$ 
or $W^d(a,b)$.
\end{Prop}
\proof
There exists a wild hypersurface bundle of degree $2$ over 
${\mathbb P}^1\times{\mathbb P}^1\times{\mathbb P}^1$ 
(see Mori-Saito \cite{mori1}). So, suppose 
$S\not\simeq{\mathbb P}^1\times{\mathbb P}^1\times
{\mathbb P}^1$. 

Suppose $\# PC(\Sigma)=5$, that is, the case (ii) in 
Theorem \ref{tmjtmj}. We use the same notation as  in Theorem 
\ref{tmjtmj}. First, we remark that 
the first, second and fourth 
primitive relations are extremal. By Corollary \ref{torickey}, 
we have $p_0+p_1+p_2-1+p_3=p_1+p_2+p_4=p_3+p_4+p_1=d+1$. 
This is impossible. 

So, we have $\# PC(\Sigma)=3$. There exists a primitive relation 
$x_{d+2}+x_{d+3}=0$ by Corollary \ref{torickey}, since 
$S$ is projective. In particular, $p=2$. 
Suppose that the types of the other 
extremal rays corresponding to 
the primitive relations $P_1$ and $P_2$ are small. 
Then, $\sigma(P_1)\cap P_2\neq\emptyset$, 
$\sigma(P_1)\cap\{x_{d+2},x_{d+3}\}\neq\emptyset$, 
$\sigma(P_2)\cap P_1\neq\emptyset$ and 
$\sigma(P_2)\cap\{x_{d+2},x_{d+3}\}\neq\emptyset$. 
This is impossible, because 
$S$ must be a ${\mathbb P}^1$-bundle over a 
toric $(d-1)$-fold with Picard number $2$.
Thus, we have the primitive 
relation $x_d+x_{d+1}=0$. It is obvious that 
the last primitive relation is 
$x_1+\cdots+x_{d-1}=\alpha x_d+\beta x_{d+2}$. 
Moreover, by Corollary \ref{torickey}, 
$\alpha$ and $\beta$ are odd numbers.
\qed

\bigskip

Let $S=W^d(a,b)$ and $D_1,\ldots,
D_{d+3}$ be the torus invariant prime divisors corresponding 
to $x_1,\ldots,x_{d+3}$, respectively. 
Put $$E\simeq\mathcal{O}_S^{\oplus (d-1)}\oplus
\mathcal{O}_S((a-1)D_1+D_{d+1})\oplus
\mathcal{O}_S((b-1)D_1+D_{d+1})\ \mbox{and}\ 
L\simeq\mathcal{O}_S(D_1).$$
We can construct a wild hypersurface bundle for these $E$ 
and $L$ similarly as in the case (I).

Let $\widetilde{\Sigma}$ be the fan corresponding to 
$Y=\mathbb{P}_S(E)$. The primitive relations of 
$\widetilde{\Sigma}$ are $\widetilde{x}_d+
\widetilde{x}_{d+1}=y_{1}$, $\widetilde{x}_{d+2}+
\widetilde{x}_{d+3}=y_{2}$,
$y_{1}+\cdots+y_{d+1}=0$ and
$$\widetilde{x}_{1}+\cdots+\widetilde{x}_{d-1}
=(2a-1)\widetilde{x}_{d}+(2b-1)\widetilde{x}_{d+2}
+(a-1)y_1+(b-1)y_2+y_{3}+\cdots+y_{d+1}$$
if $a=1$ or $b=1$, otherwise 
$$\widetilde{x}_{1}+\cdots+\widetilde{x}_{d-1}
=(2a-1)\widetilde{x}_{d}+(2b-1)\widetilde{x}_{d+2}
+(a-2)y_{1}+(b-2)y_2,$$
where $\G(\widetilde{\Sigma})=\{\widetilde{x}_1,\ldots,
\widetilde{x}_{d+3},y_1,\ldots,
y_{d+1}\}$. Let $\widetilde{D}_1,\ldots,
\widetilde{D}_{d+3},F_1,\ldots,
F_{d+1}$ be the torus invariant prime 
divisors corresponding to $\widetilde{x}_1,\ldots,
\widetilde{x}_{d+3},y_1,\ldots,
y_{d+1}$, respectively. 
Then, we have 
$\widetilde{D}_1=\cdots=\widetilde{D}_{d-1}$, 
$\widetilde{D}_{d+1}=(2a-1)\widetilde{D}_1+
\widetilde{D}_{d}$, 
$\widetilde{D}_{d+3}=(2b-1)\widetilde{D}_1+
\widetilde{D}_{d+2}$, $F_3=\cdots=F_{d+1}$ and 
$F_{d+1}=(a-1)\widetilde{D}_1+\widetilde{D}_{d}
+F_1=(b-1)\widetilde{D}_1+\widetilde{D}_{d+2}
+F_2$ in $\Pic Y$. 
Since the tautological line bundle $\xi$ for $\pi:Y\to S$ 
is $\mathcal{O}_Y(F_{r+1})$, we have 
$X\sim 2\xi+\pi^{*}L=2F_{d+1}+\widetilde{D}_1
=\widetilde{D}_{d}+\widetilde{D}_{d+1}+
2F_{1}=\widetilde{D}_{d+2}+\widetilde{D}_{d+3}+
2F_{2}$. 
Thus, for example, the smooth 
hypersurface $X$ in $Y$ defined by the equation
$$X_{d}X_{d+1}Y_1^2+X_{d+2}X_{d+3}Y_2^2
+X_1Y_3^2+\cdots+X_{d-1}Y_{d+1}^2=0$$
is a wild hypersurface 
bundle of degree $2$ over $S$, 
where $X_1,\ldots,X_{d+3},Y_1,\ldots,Y_{d+1}$ are 
the homogeneours coordinates of $Y$ corresponding to
$\widetilde{D}_1,\ldots,
\widetilde{D}_{d+3},F_1,\ldots,
F_{d+1}$, respectively. We can easily check the 
smoothness of $X$, so we leave the details for the exercise.

\section{Toric Fano varieties}\label{overfano}

\hspace{5mm} In this section, we consider the case 
where $S$ is 
a toric Fano $d$-fold. A {\em Fano} variety is 
a Gorenstein projective 
variety $S$ whose anti-canonical divisor $-K_S$ is ample. 
We can easily check 
whether a given smooth projective toric 
variety is Fano or not using the notion of primitive 
collections and relations.

\begin{Prop}[Batyrev \cite{batyrev4}, Sato \cite{sato1}]
Let $S=S_{\Sigma}$ be a smooth projective toric variety. 
$S$ is a Fano variety if and only if $\deg P>0$ for any 
primitive collection $P\in\PC(\Sigma)$.
\end{Prop}

Smooth toric Fano $d$-folds are 
classified for $d\leq 4$ (see Batyrev \cite{batyrev1}, 
\cite{batyrev4}, Oda \cite{oda2}, Sato \cite{sato1} and 
Watanabe-Watanabe \cite{watanabe1}). So, we determine 
the possibilities for these classified toric Fano 
varieties and construct wild hypersurface bundles over 
them. 

The following Proposition is easy. 

\begin{Prop}\label{divcont}
Let $f:X\to S$ be a wild hypersurface bundle over 
a toric Fano $d$-fold $S=S_{\Sigma}$ and $d\geq 3$. 
If there exists an extremal divisorial contraction 
$\varphi:S\to\overline{S}$, then 
$$S\simeq{\mathbb P}_{{\mathbb P}^{d-1}}
(\mathcal{O}_{{\mathbb P}^{d-1}}\oplus
\mathcal{O}_{{\mathbb P}^{d-1}}(2a-1))$$
for a positive integer $a$. 
\end{Prop}
\proof
By Corollary \ref{torickey}, the image of the exceptional 
divisor of $\varphi$ is a point. So, there exist 
exactly two cases 
by Bonavero's classification (see Bonavero \cite{bonavero1}): 
(a) The Picard number of $S$ is two, or (b) The Picard number 
of $S$ is three and $\#PC(\Sigma)=5$. 
However, the case (b) does not 
occur by Proposition \ref{toricpic3}. Thus, we complete the 
proof by Proposition \ref{toricpic2}.
\qed

\begin{Cor}\label{fanomain}
Let $f:X\to S$ be a wild hypersurface bundle over 
a toric Fano $d$-fold $S=S_{\Sigma}$ and $d\geq 3$. 
Then, one of the following holds$:$
\begin{enumerate}
\item $S\simeq{\mathbb P}^d$. 
\item $S\simeq({\mathbb P}^1)^{d}$.
\item $S\simeq{\mathbb P}_{{\mathbb P}^{d-1}}
(\mathcal{O}_{{\mathbb P}^{d-1}}\oplus
\mathcal{O}_{{\mathbb P}^{d-1}}(2a-1))$
for a positive integer $a$.
\item Every extremal contraction of $S$ is 
either a $\mathbb{P}^1$-bundle structure or 
a small contraction, and there exists at least 
one small contraction.
\end{enumerate}
\end{Cor}
\proof
See Mori-Saito \cite{mori1} for the cases (i) and (ii), and 
see the case (I) in Section \ref{pic2or3} for the case (iii). 
So, suppose $S$ is not one of them. 
For the case (i) in Corollary \ref{torickey}, 
we have $n\geq 2$ 
by Proposition \ref{divcont}. For the case (ii) 
in Corollary \ref{torickey}, 
we have $n=0$ and the associated extremal contraction is 
a $\mathbb{P}^1$-bundle structure, 
since $S$ is a Fano variety.
\qed

\bigskip

(I) $\dim S=2$.

\bigskip

There exist exactly five toric del Pezzo surfaces
$${\mathbb P}^2,\ {\mathbb P}^1\times{\mathbb P}^1,
\ {\mathbb P}_{{\mathbb P}^1}\left(
\mathcal{O}_{{\mathbb P}^1}
\oplus\mathcal{O}_{{\mathbb P}^1}(1)\right),
\ S_6\ \mbox{and}\ S_7,$$
where $S_6$ and $S_7$ are the del Pezzo surfaces 
of degree $6$ and $7$, respectively. For any toric 
del Pezzo surface $S$, there exists a wild 
hypersurface bundle over $S$. In fact, for $S_6$ and $S_7$, 
we can construct wild hypersurface bundles similarly as 
in Section \ref{pic2or3}. We omit the precise calculation 
for these constructions, and use the same 
notation as in Section \ref{pic2or3}. 

\begin{Ex}
{\rm
Let $S=S_{\Sigma}$ be the del Pezzo surface 
$S_7$ of degree $7$. The primitive relations are
$x_1+x_2=x_3$, $x_1+x_5=0$, $x_2+x_4=x_5$, $x_3+x_4=0$ 
and $x_3+x_5=x_2$. We have $p=2$. 
Put 
$$E=\mathcal{O}_S\oplus\mathcal{O}_S(D_3)\oplus
\mathcal{O}_S(D_5)\ \mbox{and}
\ L=\mathcal{O}_S(D_2).$$ The 
primitive relations of $\widetilde{\Sigma}$ are 
$\widetilde{x}_1+\widetilde{x}_2=\widetilde{x}_3+y_2+y_3$, 
$\widetilde{x}_1+\widetilde{x}_5=y_2$, 
$\widetilde{x}_2+\widetilde{x}_4=\widetilde{x}_5+y_1+y_3$, 
$\widetilde{x}_3+\widetilde{x}_4=y_1$, 
$\widetilde{x}_3+\widetilde{x}_5=\widetilde{x}_2+y_1+y_2$ 
and $y_1+y_2+y_3=0$. 
The hypersurface $X$ in $Y=\mathbb{P}_S(E)$ 
defined by the equation
$$X_3X_4Y_1^2+X_1X_5Y_2^2+X_2Y_3^2=0$$
is a wild hypersurface bundle of degree $2$ over $S$.
}
\end{Ex}

\begin{Ex}
{\rm
Let $S=S_{\Sigma}$ be the del Pezzo surface 
$S_6$ of degree $6$. The primitive relations are 
$x_1+x_5=0$, $x_3+x_4=0$, $x_2+x_6=0$, 
$x_3+x_6=x_1$, $x_3+x_5=x_2$, $x_1+x_2=x_3$, 
$x_5+x_6=x_4$, $x_2+x_4=x_5$ and $x_1+x_4=x_6$. 
We have $p=2$. 
Put 
$$E=\mathcal{O}_S\oplus\mathcal{O}_S(D_5-D_6)\oplus
\mathcal{O}_S(-D_2+D_4)\ \mbox{and}
\ L=\mathcal{O}_S(D_2+D_3).$$
The primitive relations of $\widetilde{\Sigma}$ are 
$\widetilde{x}_1+\widetilde{x}_5=y_1$, 
$\widetilde{x}_3+\widetilde{x}_4=y_2$, 
$\widetilde{x}_2+\widetilde{x}_6=y_3$, 
$\widetilde{x}_3+\widetilde{x}_6=\widetilde{x}_1+y_2+y_3$, 
$\widetilde{x}_3+\widetilde{x}_5=\widetilde{x}_2+y_1+y_2$, 
$\widetilde{x}_1+\widetilde{x}_2=\widetilde{x}_3+y_1+y_3$, 
$\widetilde{x}_5+\widetilde{x}_6=\widetilde{x}_4+y_1+y_3$, 
$\widetilde{x}_2+\widetilde{x}_4=\widetilde{x}_5+y_2+y_3$, 
$\widetilde{x}_1+\widetilde{x}_4=\widetilde{x}_6+y_1+y_2$ 
and $y_1+y_2+y_3=0$. 
The hypersurface $X$ in $Y=\mathbb{P}_S(E)$ 
defined by the equation
$$X_1X_5Y_1^2+X_3X_4Y_2^2+X_2X_6Y_3^2=0$$
is a wild hypersurface bundle of degree $2$ over $S$.
}
\end{Ex}

(II) $\dim S=3$.

\bigskip

There does not exist a small contraction from any 
smooth toric Fano $3$-fold. Therefore, if there 
exists a wild hypersurface bundle over $S$, then 
$S$ is isomorphic to one of the following 
by Corollary \ref{fanomain}:
$${\mathbb P}^3,\ {\mathbb P}^1\times{\mathbb P}^1
\times{\mathbb P}^1\ \mbox{and}
\ {\mathbb P}_{{\mathbb P}^2}\left(
\mathcal{O}_{{\mathbb P}^2}
\oplus\mathcal{O}_{{\mathbb P}^2}(1)\right).$$

\bigskip

(III) $\dim S=4$.

\bigskip

There exists a wild hypersurface bundle over $S$, 
if $S$ is isomorphic to one of the following:
$${\mathbb P}^4,\ {\mathbb P}^1\times{\mathbb P}^1
\times{\mathbb P}^1\times{\mathbb P}^1,
\ {\mathbb P}_{{\mathbb P}^3}\left(
\mathcal{O}_{{\mathbb P}^3}
\oplus\mathcal{O}_{{\mathbb P}^3}(1)\right)\ \mbox{and}
\ {\mathbb P}_{{\mathbb P}^3}\left(
\mathcal{O}_{{\mathbb P}^3}
\oplus\mathcal{O}_{{\mathbb P}^3}(3)\right).$$
So, suppose $S$ is not one of them, that is, the case 
(iv) in Corollary \ref{fanomain}. By the classification of 
smooth toric Fano $4$-folds, there exist exactly four 
possibilities:
\begin{enumerate}
\item $S\simeq W^4(1,1),$
\item $S$ is the toric Fano $4$-fold of type $M_1$ 
(see Batyrev \cite{batyrev4} and Sato \cite{sato1}),
\item $S$ is the $4$-dimensional pseudo del Pezzo 
variety $\widetilde{V}^4$ (see Ewald \cite{ewald1}) and
\item $S$ is the $4$-dimensional del Pezzo 
variety $V^4$ (see Klyachko-Voskresenskij \cite{voskre1}).
\end{enumerate}
The first case is studied in Section \ref{pic2or3}, and 
we can construct wild hypersurface bundles for the other 
cases similarly as 
in Section \ref{pic2or3}. We omit the precise calculation 
for these constructions, and use the same notation 
as in Section \ref{pic2or3}. 

\begin{Ex}
{\rm
Let $S=S_{\Sigma}$ be the toric Fano $4$-fold of type $M_1$. 
The primitive relations are 
$x_1+x_8=0$, $x_4+x_5=0$, $x_6+x_7=0$, 
$x_1+x_2+x_3=x_4+x_6$, $x_4+x_6+x_8=x_2+x_3$, 
$x_2+x_3+x_5=x_6+x_8$ and $x_2+x_3+x_7=x_4+x_8$. 
We have $p=2$. 
Put 
$$E=\mathcal{O}_S\oplus\mathcal{O}_S\oplus
\mathcal{O}_S(D_8)\oplus\mathcal{O}_S(D_4)\oplus
\mathcal{O}_S(D_6)\ \mbox{and}\ L=\mathcal{O}_S(D_3).$$ 
The primitive relations of $\widetilde{\Sigma}$ are 
$\widetilde{x}_1+\widetilde{x}_8=y_1$, 
$\widetilde{x}_4+\widetilde{x}_5=y_2$, 
$\widetilde{x}_6+\widetilde{x}_7=y_3$, 
$\widetilde{x}_1+\widetilde{x}_2+\widetilde{x}_3=
\widetilde{x}_4+\widetilde{x}_6+y_1+y_4+y_5$, 
$\widetilde{x}_4+\widetilde{x}_6+\widetilde{x}_8=
\widetilde{x}_2+\widetilde{x}_3+y_1+y_2+y_3$, 
$\widetilde{x}_2+\widetilde{x}_3+\widetilde{x}_5=
\widetilde{x}_6+\widetilde{x}_8+y_2+y_4+y_5$, 
$\widetilde{x}_2+\widetilde{x}_3+\widetilde{x}_7=
\widetilde{x}_4+\widetilde{x}_8+y_3+y_4+y_5$ 
and $y_1+y_2+y_3+y_4+y_5=0$. 
The hypersurface $X$ in $Y=\mathbb{P}_S(E)$ 
defined by the equation
$$X_1X_8Y_1^2+X_4X_5Y_2^2+X_6X_7Y_3^2+
X_2Y_4^2+X_3Y_5^2=0$$
is a wild hypersurface bundle of degree $2$ over $S$.
} 
\end{Ex}

\begin{Ex}
{\rm
Let $S=S_{\Sigma}$ be the $4$-dimensional pseudo del Pezzo 
variety $\widetilde{V}^4$. The primitive relations are 
$x_4+x_9=0$, $x_1+x_5=0$, $x_2+x_6=0$, $x_3+x_7=0$, 
$x_1+x_2+x_9=x_7+x_8$, $x_1+x_3+x_9=x_6+x_8$, 
$x_2+x_3+x_9=x_5+x_8$, $x_1+x_2+x_3=x_4+x_8$, 
$x_4+x_5+x_8=x_2+x_3$, $x_4+x_6+x_8=x_1+x_3$, 
$x_4+x_7+x_8=x_1+x_2$, $x_5+x_6+x_8=x_3+x_9$, 
$x_5+x_7+x_8=x_2+x_9$ and $x_6+x_7+x_8=x_1+x_9$. 
We have $p=2$. 
Put 
$$E=\mathcal{O}_S\oplus\mathcal{O}_S(D_1)\oplus
\mathcal{O}_S(D_2)\oplus\mathcal{O}_S(D_3)\oplus
\mathcal{O}_S(D_9)\ \mbox{and}\ L=\mathcal{O}_S(D_8).$$
The primitive relations of $\widetilde{\Sigma}$ are 
$\widetilde{x}_4+\widetilde{x}_9=y_4$, 
$\widetilde{x}_1+\widetilde{x}_5=y_1$, 
$\widetilde{x}_2+\widetilde{x}_6=y_2$, 
$\widetilde{x}_3+\widetilde{x}_7=y_3$, 
$\widetilde{x}_1+\widetilde{x}_2+\widetilde{x}_9=
\widetilde{x}_7+\widetilde{x}_8+y_1+y_2+y_4$, 
$\widetilde{x}_1+\widetilde{x}_3+\widetilde{x}_9=
\widetilde{x}_6+\widetilde{x}_8+y_1+y_3+y_4$, 
$\widetilde{x}_2+\widetilde{x}_3+\widetilde{x}_9=
\widetilde{x}_5+\widetilde{x}_8+y_2+y_3+y_4$, 
$\widetilde{x}_1+\widetilde{x}_2+\widetilde{x}_3=
\widetilde{x}_4+\widetilde{x}_8+y_1+y_2+y_3$, 
$\widetilde{x}_4+\widetilde{x}_5+\widetilde{x}_8=
\widetilde{x}_2+\widetilde{x}_3+y_1+y_4+y_5$, 
$\widetilde{x}_4+\widetilde{x}_6+\widetilde{x}_8=
\widetilde{x}_1+\widetilde{x}_3+y_2+y_4+y_5$, 
$\widetilde{x}_4+\widetilde{x}_7+\widetilde{x}_8=
\widetilde{x}_1+\widetilde{x}_2+y_3+y_4+y_5$, 
$\widetilde{x}_5+\widetilde{x}_6+\widetilde{x}_8=
\widetilde{x}_3+\widetilde{x}_9+y_1+y_2+y_5$, 
$\widetilde{x}_5+\widetilde{x}_7+\widetilde{x}_8=
\widetilde{x}_2+\widetilde{x}_9+y_1+y_3+y_5$, 
$\widetilde{x}_6+\widetilde{x}_7+\widetilde{x}_8=
\widetilde{x}_1+\widetilde{x}_9+y_2+y_3+y_5$ 
and $y_1+y_2+y_3+y_4+y_5=0$. 
The hypersurface $X$ in $Y=\mathbb{P}_S(E)$ 
defined by the equation
$$X_1X_5Y_1^2+X_2X_6Y_2^2+X_3X_7Y_3^2+
X_4X_9Y_4^2+X_8Y_5^2=0$$
is a wild hypersurface bundle of degree $2$ over $S$.
}
\end{Ex}

\begin{Ex}
{\rm
Let $S=S_{\Sigma}$ be the $4$-dimensional del Pezzo 
variety $V^4$. The primitive relations are 
$x_4+x_{10}=0$, $x_1+x_5=0$, $x_2+x_6=0$, $x_3+x_7=0$, 
$x_8+x_9=0$, $x_1+x_2+x_{10}=x_7+x_8$, 
$x_1+x_3+x_{10}=x_6+x_8$, $x_2+x_3+x_{10}=x_5+x_8$, 
$x_1+x_2+x_3=x_4+x_8$, $x_1+x_9+x_{10}=x_6+x_7$, 
$x_2+x_9+x_{10}=x_5+x_7$, $x_3+x_9+x_{10}=x_5+x_6$, 
$x_1+x_2+x_9=x_4+x_7$, $x_1+x_3+x_9=x_4+x_6$, 
$x_2+x_3+x_9=x_4+x_5$, $x_4+x_5+x_6=x_3+x_9$, 
$x_4+x_5+x_7=x_2+x_9$, $x_4+x_6+x_7=x_1+x_9$, 
$x_5+x_6+x_7=x_9+x_{10}$, $x_4+x_5+x_8=x_2+x_3$, 
$x_4+x_6+x_8=x_1+x_3$, $x_4+x_7+x_8=x_1+x_2$, 
$x_5+x_6+x_8=x_3+x_{10}$, $x_5+x_7+x_8=x_2+x_{10}$ 
and $x_6+x_7+x_8=x_1+x_{10}$. 
We have $p=2$. 
Put 
$$E=\mathcal{O}_S\oplus\mathcal{O}_S(D_1-D_9)\oplus
\mathcal{O}_S(D_2-D_9)\oplus\mathcal{O}_S(D_3-D_9)\oplus
\mathcal{O}_S(D_{10}-D_9)\ \mbox{and}$$
$$L=\mathcal{O}_S(D_8+D_9).$$
The primitive relations of $\widetilde{\Sigma}$ are 
$\widetilde{x}_4+\widetilde{x}_{10}=y_4$, 
$\widetilde{x}_1+\widetilde{x}_5=y_1$, 
$\widetilde{x}_2+\widetilde{x}_6=y_2$, 
$\widetilde{x}_3+\widetilde{x}_7=y_3$, 
$\widetilde{x}_8+\widetilde{x}_9=y_5$, 
$\widetilde{x}_1+\widetilde{x}_2+\widetilde{x}_{10}=
\widetilde{x}_7+\widetilde{x}_8+y_1+y_2+y_4$, 
$\widetilde{x}_1+\widetilde{x}_3+\widetilde{x}_{10}=
\widetilde{x}_6+\widetilde{x}_8+y_1+y_3+y_4$, 
$\widetilde{x}_2+\widetilde{x}_3+\widetilde{x}_{10}=
\widetilde{x}_5+\widetilde{x}_8+y_2+y_3+y_4$, 
$\widetilde{x}_1+\widetilde{x}_2+\widetilde{x}_3=
\widetilde{x}_4+\widetilde{x}_8+y_1+y_2+y_3$, 
$\widetilde{x}_1+\widetilde{x}_9+\widetilde{x}_{10}=
\widetilde{x}_6+\widetilde{x}_7+y_1+y_4+y_5$, 
$\widetilde{x}_2+\widetilde{x}_9+\widetilde{x}_{10}=
\widetilde{x}_5+\widetilde{x}_7+y_2+y_4+y_5$, 
$\widetilde{x}_3+\widetilde{x}_9+\widetilde{x}_{10}=
\widetilde{x}_5+\widetilde{x}_6+y_3+y_4+y_5$, 
$\widetilde{x}_1+\widetilde{x}_2+\widetilde{x}_9=
\widetilde{x}_4+\widetilde{x}_7+y_1+y_2+y_5$, 
$\widetilde{x}_1+\widetilde{x}_3+\widetilde{x}_9=
\widetilde{x}_4+\widetilde{x}_6+y_1+y_3+y_5$, 
$\widetilde{x}_2+\widetilde{x}_3+\widetilde{x}_9=
\widetilde{x}_4+\widetilde{x}_5+y_2+y_3+y_5$, 
$\widetilde{x}_4+\widetilde{x}_5+\widetilde{x}_6=
\widetilde{x}_3+\widetilde{x}_9+y_1+y_2+y_4$, 
$\widetilde{x}_4+\widetilde{x}_5+\widetilde{x}_7=
\widetilde{x}_2+\widetilde{x}_9+y_1+y_3+y_4$, 
$\widetilde{x}_4+\widetilde{x}_6+\widetilde{x}_7=
\widetilde{x}_1+\widetilde{x}_9+y_2+y_3+y_4$, 
$\widetilde{x}_5+\widetilde{x}_6+\widetilde{x}_7=
\widetilde{x}_9+\widetilde{x}_{10}+y_1+y_2+y_3$, 
$\widetilde{x}_4+\widetilde{x}_5+\widetilde{x}_8=
\widetilde{x}_2+\widetilde{x}_3+y_1+y_4+y_5$, 
$\widetilde{x}_4+\widetilde{x}_6+\widetilde{x}_8=
\widetilde{x}_1+\widetilde{x}_3+y_2+y_4+y_5$, 
$\widetilde{x}_4+\widetilde{x}_7+\widetilde{x}_8=
\widetilde{x}_1+\widetilde{x}_2+y_3+y_4+y_5$, 
$\widetilde{x}_5+\widetilde{x}_6+\widetilde{x}_8=
\widetilde{x}_3+\widetilde{x}_{10}+y_1+y_2+y_5$, 
$\widetilde{x}_5+\widetilde{x}_7+\widetilde{x}_8=
\widetilde{x}_2+\widetilde{x}_{10}+y_1+y_3+y_5$, 
$\widetilde{x}_6+\widetilde{x}_7+\widetilde{x}_8=
\widetilde{x}_1+\widetilde{x}_{10}+y_2+y_3+y_5$ 
and $y_1+y_2+y_3+y_4+y_5=0$. 
The hypersurface $X$ in $Y=\mathbb{P}_S(E)$ 
defined by the equation
$$X_1X_5Y_1^2+X_2X_6Y_2^2+X_3X_7Y_3^2+
X_4X_{10}Y_4^2+X_8X_9Y_5^2=0$$
is a wild hypersurface bundle of degree $2$ over $S$.
}
\end{Ex}


\bigskip

\begin{flushleft}
\begin{sc}
Department of Mathematics \\
Tokyo Institute of Technology \\
2-12-1 Oh-okayama \\
Meguro-ku \\
Tokyo 152-8551 \\
Japan
\end{sc}

\medskip
{\it E-mail address}: $\mathtt{hirosato@math.titech.ac.jp}$
\end{flushleft}

\end{document}